\documentclass[11pt,a4paper]{article}
\usepackage{amsmath,amsthm,enumerate}
\usepackage{amssymb}
\usepackage{amsfonts}
\usepackage{graphicx}
\usepackage{color}
\usepackage{csquotes}
\usepackage{tcolorbox}

\newtheorem{Theorem}{Theorem}[section]

\newtheorem{Proposition}{Proposition}[section]

\newtheorem{Lemma}{Lemma}[section]

\def\bc{\begin{center}}
\def\ec{\end{center}}

\def\s2c{\vskip 2cm}

\def\R{\mathbb{R}}

\def\N{\mathbb{N}}

\def\rar{\rightarrow}

\def\ve{\varepsilon}
\numberwithin{equation}{section}

\def\bt{\begin{Theorem}}
\def\et{\end{Theorem}}
\def\bl{\begin{Lemma}}
\def\el{\end{Lemma}}
\def\bcor{\begin{Corollary}}
\def\ecor{\end{Corollary}}

\begin{document}

\title{Simple Bilevel Programming and Extensions Part-II: Algorithms}
\author{Stephan Dempe\footnote{Faculty of Mathematics and Computer Science, TU Bergakademie Freiberg, Germany, e-mail: dempe@tu-freiberg.de. Work of this author has been supported by Deutsche Forschungsgemeinschaft}, Nguyen Dinh\footnote{Department of Mathematics, International University, Vietnam National University-Ho Chi Minh City,
		Ho Chi Minh City, Vietnam, e-mail: ndinh02@gmail.com.  \ \ Work of this author was supported by the project 101.01-2015.27
 \emph{Generalizations of Farkas lemma with applications to optimization}    from the National Foundation for Science \& Technology Development (NAFOSTED),
Vietnam.
}, Joydeep Dutta\footnote{Department of Economic Sciences, Indian
Institute of Technology, Kanpur, India, e-mail: jdutta@iitk.ac.in} and Tanushree Pandit\footnote{Department of Mathematics and Statistics, Indian
Institute of Technology, Kanpur, India, e-mail: tpandit@iitk.ac.in}}
\date{}
\maketitle

\begin{abstract}
This article continues our study on simple bilevel and simple MPEC problems. In this article we focus on developing algorithms. We show how using the idea of a gap function one can represent a simple MPEC as a simple bilevel problem with non-smooth data. This motivates us first to develop an algorithm for a simple bilevel problem with non-smooth data and modify the scheme for the same to develop an algorithm for the simple MPEC problem. We also discuss how the simple bilevel formulation of a simple MPEC can help us in formulating a stopping criteria for the simple MPEC problem.	

\end{abstract}

\\
\medskip


\section{Introduction}
In this article we continue with our study of simple bilevel programming problem (SBP for short) and the simple MPEC problem (SMPEC for short), which was initiated in part-I \cite{DDDP}. In this article which we mark as part-II, our focus will be on algorithms. Our main aim would be to develop an algorithm for the (SMPEC) problem. Our motivation comes from the last section of part-I \cite{DDDP}. There in we presented a schematic algorithm for the (SMPEC) problem, whose convergence analysis depended on sequential optimality conditions. We had called this algorithm in \cite{DDDP} schematic since the dual gap function of the variational inequality played a central role there. The dual gap function is often difficult to compute and thus an attempt is made in this paper to develop an algorithm, where we do not have to involve the dual gap function. However, as we will observe, avoiding the dual gap function has its own cost, in the form of an additional assumption. In fact for convergence analysis we need the assumption of monotone plus on the vector function associated with the variational inequality at the lower level of the (SMPEC) problem. We want also to mention, that from the perspective of theoretical discussions the dual gap function will continue to play an important role even in this paper. In fact, as we will discuss at the end of the paper that if suitable approaches to compute the dual gap function can be developed more cheaply, then it may play an important role in the building of implementable algorithms as well.\\
In this paper we develop algorithms for both the (SBP) and (SMPEC) problem. For the (SBP) problem we will consider both the upper and lower level objective functions to be non-smooth. As we will discuss below that our motivation to develop an algorithm for such a problem comes from the reformulation of the (SMPEC) as an (SBP) using the dual gap function. Our approach to solve the (SBP) problem with non-smooth data will be a mixture of penalization and proximal point techniques and the scheme would be suitably modified to fit the (SMPEC) problem.\\
Let us just recall that an (SMPEC) problem is given as follows
\begin{eqnarray*}
&&\min f(x)\\
&&\mbox{subject to } x\in \mbox{sol(VI}(F,C)),
\end{eqnarray*}
where $ f : \R^n \rightarrow \R $ is a convex function, which need not be differentiable throughout and $ \mbox{VI}(F,C) $ represents as usual the variational inequality described by a vector function $ F : \R^n \rightarrow \R^n $ and a closed convex set $ C $, where in we seek to find an $ x\in C $ such that
\begin{eqnarray*}
\langle F(x), y-x \rangle \geq 0 \quad \forall y \in C
\end{eqnarray*}
and $\mbox{sol(VI}(F,C)) $ denotes the set of its solutions. Additionally we assume here that $ F : \R^n \rightarrow \R^n $ is continuous and monotone, {\it i.e.}
\begin{eqnarray*}
\langle F(y)-F(x), y-x \rangle \geq 0 \quad \forall x,y \in \R^n.
\end{eqnarray*}
The assumption that $ F $ is monotone guarantees that $ \mbox{sol(VI}(F,C)) $ is convex. We would again recall from part-I \cite{DDDP} that the dual gap function $ g_D $ associated with $ \mbox{VI}(F,C) $ is given as
\begin{eqnarray*}
g_D(x)= \sup\limits_{y \in C} \langle F(y), x-y \rangle.
\end{eqnarray*}
We will work under the assumption that $\mbox{sol(VI}(F,C)) \neq \emptyset $. Under this assumption we know from section 3.2 of part-I \cite{DDDP} that
\begin{eqnarray*}
\mbox{sol(VI}(F,C)) = \mbox{argmin}_{C} g_D.
\end{eqnarray*}
Thus the (SMPEC) problem can be equivalently written as an (SBP) problem as follows,
\begin{eqnarray*}
&&\min f(x)\\
&&\mbox{subject to } x\in \mbox{argmin}_{C} g_D.
\end{eqnarray*}
Let us call the above problem (SBP-1). Observe that (SBP-1) is an (SBP) problem where both upper and lower level objectives are non-smooth. Note that $ g_D$ in general is not smooth. This motivates us to develop an algorithm for the (SBP) problem with non-smooth data. Though we have discussed the main tools from convex analysis used in our discussion in Part-I \cite{DDDP}, we recall here the notion of $\varepsilon$-subdifferential, $ \varepsilon $-normal set and the subdifferential sum rule for the sake of completeness.\\

The $\varepsilon $-{\it subdifferential} of $f$, $\partial
_{\varepsilon }f: \R^n \rightrightarrows \R^n$,  is defined as
\begin{equation*}
\partial _{\varepsilon }{f}(x)=\{v\in \R^n\; :\;  f(y)\geq
f(x)+\langle v,y-x\rangle - \varepsilon ,\,\forall \,y\in \R^n \}.
\end{equation*}
It is important to note that if $f$ is a proper and lower-semicontinuous function then $\partial_{\varepsilon } f(x) \neq \emptyset$  for any $ x \in {\rm dom} f$ and any $\varepsilon >0$. The $\varepsilon
$-subdifferential is a key tool to analyze the notion of an $\varepsilon$-minimizer of a proper, lsc and  convex function. Given $\varepsilon >0$, a vector $\bar{x} \in \mathbb{R}^n$ is called an $\varepsilon$-minimizer if $f(\bar{x}) \leq f(x) + \varepsilon$ for all $x \in \mathbb{R}^n$. The set of all $\varepsilon$-minimizers of $f$ is denoted as $\varepsilon\mbox{-argmin}f$. In fact $\bar{x}$ is an $\varepsilon$-minimizer of a proper, lsc, convex function $f$ if and only if $0 \in \partial_{\varepsilon} f(\bar{x})$. This notion of an $\varepsilon$-minimizer can be extended to cover the constrained case. Of course a vector $\bar{x} \in \mathbb{R}^n$ is an $\varepsilon$-minimizer of a convex function $f$ over a convex set $C$, if $f(\bar{x}) \leq f(x)+\varepsilon$ for all $x \in C$. In order to develop  necessary and sufficient optimality conditions for an $\varepsilon$-minimizer of $f$ over $C$, we need to know the sum rule for $\varepsilon$-subdifferentials and also the notion of an $\varepsilon$-normal set. The $\varepsilon$-normal set of a convex set $C$ at any $\bar{x} \in C$ is given as
\begin{eqnarray*}
N_C^{\varepsilon}(\bar{x})= \{v \in \mathbb{R}^n: \langle v, x- \bar{x}\rangle \leq \varepsilon, \quad \forall x \in C\}.
\end{eqnarray*}
In fact $\partial_{\varepsilon} \delta_C(\bar{x}) = N _C ^{\varepsilon} (\bar{x})$. Let $f$ and $ g $ be any convex functions and let $\varepsilon > 0$ be given. Then
\begin{eqnarray*}
\partial_{\varepsilon}(f + g)(x) \subset \bigcup \limits_{\substack{\varepsilon= \varepsilon_1 + \varepsilon_2\\ \varepsilon_1 >0, \varepsilon_2 >0}} (\partial_{\varepsilon_1} f(x) + \partial_{\varepsilon_2} g(x)).
\end{eqnarray*}
For more details on the sum rule see Hiriart-Urruty and Lemarechal \cite{HU-Lemarechal}. A vector $\bar{x} \in \mathbb{R}^n$ is an $\varepsilon$-minimizer of a convex function $f $ over a convex set $C$ if and only if there exists $\varepsilon_1 >0$ and $\varepsilon_2 >0$ such that $\varepsilon_1 + \varepsilon_2 = \varepsilon$ and
\begin{eqnarray}\label{2.2}
0 \in \partial_{\varepsilon_1}f(\bar{x}) + N_C^{\varepsilon_2}(\bar{x}).
\end{eqnarray}
Furthermore,  if $f$ is differentiable at $\bar{x}$, then $\nabla f (\bar{x}) \in \partial _{\varepsilon} f(\bar{x})$ for any $\epsilon>0$, but $\partial _{\varepsilon} f(\bar{x})$ need not be a singleton. This property along with (\ref{2.2}) will play a fundamental role in designing our algorithm for the simple bilevel problem with non-smooth data.\\
In our study of the convergence analysis for the algorithm that is developed here for a simple bilevel programming problem, we shall need the notion of a recession cone or asymptotic cone of a convex set and recession function of a convex function.\\
Given a nonempty convex set $A \subset \mathbb{R}^n$, the recession cone or asymptotic cone is given by
\begin{eqnarray}
A^{\infty}:= \{y \in \mathbb{R}^n : \forall x \in A, \forall \lambda \geq 0: x+\lambda y \in A\}.
\end{eqnarray}
The recession function of a convex function $f: \mathbb{R}^n \rightarrow \mathbb{R}$ is the function $f^{\infty}: \mathbb{R}^n \rightarrow (-\infty,+\infty]$ whose epigraph is the recession cone of the epigraph of $f$.
Also  we are using the following well-known results (for example see \cite{rocw}) regarding recession cone and recession function.
\begin{enumerate}
\item $(A \cap B)^{\infty}= A^{\infty} \cap B^{\infty}$, for any two convex sets $A,B \subset \mathbb{R}^n$.
\item A convex set $A \subset \mathbb{R}^n$ is bounded if and only if $A^{\infty} = \{0\}$.
\item For any convex function $f: \mathbb{R}^n \rightarrow \mathbb{R}$ and $M \in \mathbb{R}$,
\begin{eqnarray*}
\{x \in \mathbb{R}^n : f(x) \leq M\}^{\infty} = \{x \in \mathbb{R}^n : f^{\infty}(x) \leq 0\}.
\end{eqnarray*}
\end{enumerate}
To develop the algorithm for the (SMPEC) problem, we have assumed that the function $F$ in the lower level problem $VI(F,C)$ is monotone plus.\\
A monotone                                                                                                                                                                                                                                                                                                                                                                                                                                                                                                                                                                                                                                                                                                                                                                                                                                                                                                                                                                                                                                                                                                                                                                                                                                     function $F: \mathbb{R}^n \rightarrow \mathbb{R}^n$ is called {\it monotone plus}  if for any $x,y \in \mathbb{R}^n$,
\begin{eqnarray*}
\langle F(x)-F(y),x-y \rangle =0 \Rightarrow F(x)=F(y).
\end{eqnarray*}
The notion of a distance function will play a major role in the convergence analysis of algorithms for the (SBP) problem and the (SMPEC) problem. Given a set $ C \subset \mathbb{R}^n $ and $ x \in \mathbb{R}^n $, the distance of the point $ x $ from the set $ C $ is given as
\begin{eqnarray*}
d(x,C)= \inf\limits_{y \in C} \| y-x \|.
\end{eqnarray*}
If $ C $ is closed and convex, then the distance function is a convex function and in fact there is a unique minimizer of the problem $ \inf\limits_{y \in C} \| y-x \| $. The unique minimizer $ \hat{y} \in C $(say) is also known as the projection of $ x $ on $ C $ if $ x \notin C $. We then write $ \hat{y} = P_C(x) $ and then $ d(x,C)= \| x- P_C(x) \| $.\\

\noindent Our treatment of algorithms for (SBP) and (SMPEC) is motivated from two different sources. The first source is Facchenei et al.\cite{fpang2011}. In this paper they develop an algorithm for variational inequality constrained hemivariational inequality problem, which is in short VI-C HVI$(F,C,H,\phi)$. In this problem $ F $ and $ H $ are Lipschitz continuous maps and $ \phi $ is a smooth and Lipschitz continuous, convex function and $ C $ is a closed convex set. Let us recall from \cite{fpang2011}, that in the problem VI-C HVI$(F,C,H,\phi)$ we are seeking to find an $ x\in \mbox{Sol(VI}(F,C)) $ such that
\begin{eqnarray*}
\langle H(x), y-x \rangle + \phi(y) - \phi(x) \geq 0, \quad \forall y \in \mbox{Sol(VI}(F,C)).
\end{eqnarray*}
Observe that VI-C HVI$(F,C,H,\phi)$ reduces to simple bilevel problem if $ H \equiv 0 $. In fact for the convergence analysis of our numerical scheme for SMPEC we take several ideas from the convergence analysis of Facchenei and Pang \cite{fpang2011} for the case when $ \phi $ is not a smooth convex function and need not be globally Lipschitz and $ F $ need not be Lipschitz too. The most striking fact is that when we consider the problem (SBP) where both upper and lower level objective functions are non-smooth, we can still combine the approaches of Facchenei et al. \cite{fpang2011} with the approach of Cabot \cite{Cabot}, where a simple bilevel problem was considered with unconstrained lower level problem, to develop a convergence analysis of the numerical scheme for the (SBP) problem; that we present here. We would like to recall that the lower level problem in (SBP) is constrained and the problem has non-smooth data. To the best of our knowledge Solodov \cite{Solodov1} was the first to develop an algorithmic scheme for (SBP) with smooth data, where the upper and lower level objective functions were assumed to have Lipschitz gradients. Our approach is motivated by that of Solodov \cite{Solodov1}, though our algorithm differs in a fundamental way with that of Solodov \cite{Solodov1}. We will discuss these differences in section 2.\\

This article is organised as follows. In section 2 we present an algorithm and its convergence analysis for the (SBP) problem where both the upper and lower objective functions are not necessarily differentiable and the lower level problem is a constrained convex optimization problem, with a feasible set which is a closed convex and proper subset of $ \R^n $. As we are aware that the algorithm for the non-smooth (SBP) is developed since the (SMPEC) problem that we study here can be modelled as a non-smooth (SBP) as we have seen. In section 3, we again focus on the (SBP) reformulation of the (SMPEC) problem. We try to see to what extent one can use the (SBP) reformulation to study the (SMPEC) problem. We argue that this reformulation can lead us to a stopping criteria for any algorithm for the (SMPEC) problem, once we have been able to formulate some approximate version of the Lagrange multiplier rule for the (SBP) reformulation. In section 4 we present an algorithm for the (SMPEC) problem along with its convergence analysis. We shall see that this algorithmic scheme for the (SMPEC) problem has been obtained by tweaking the scheme we have developed for the (SBP) in section 2.\\

We would like to admit that we have not carried out numerical experiments in this article. One of the primary reasons being that it might be difficult to compute non-smooth objects like the $\ve$-subdifferential or the $\ve$-normal set to a convex set. It is apparent that implementing these algorithms for a general convex programming problem may not be easy unless some special structure is considered. As we were revising this paper some progress in that direction was made by Pandit, Dutta and Rao, which we will discuss in the conclusion of this article. We would like the reader to view this paper as an attempt to develop algorithms for (SBP) and (SMPEC), where we can dispense with the Lipschitz gradient assumption and carry out the analysis in a very general framework.

\section{Algorithm for non-smooth simple bilevel problem}
In this section we will be concerned with developing an algorithm for the (SBP) problem given as
\begin{eqnarray*}
\min f(x),
\quad \mbox{subject to} \quad x \in S_0,
\end{eqnarray*}
where $S_0= \mbox{argmin}\{g(x): x \in C\}$, $ f, g : \R^n \rightarrow \R $ are convex functions and $ C $ is a closed convex set.\\
In \cite{Solodov1} Solodov considers $f$ and $g$ to be smooth convex, Lipschitz continuous functions over $C$, where $C$ is a closed, convex set. The approach taken by Solodov \cite{Solodov1} is a method of penalization of the form $g+ \ve_n f$, where $\ve_n \downarrow 0$ as $ n \rightarrow \infty $ is a penalty parameter and solves this penalized problem over $C$ by using a projected gradient technique combined with an Armijo type of decrease criteria. Solodov \cite{Solodov2} later on also developed a bundle method 
when $f$ and $g$ are nonsmooth and convex and $C=\mathbb{R}^n$. In our setting we shall consider $f$ and $g$ to be convex and nonsmooth but $C$ is a non-empty, proper, convex subset of $\mathbb{R}^n$, which may be unbounded. Our approach is to use Solodov's penalization technique in the nonsmooth setting. We should also note that, in our approach we shall use an inexact proximal point approach to solve the penalized problem in contrast to the projected gradient approach of Solodov \cite{Solodov1}. Further we shall not use any decrease criteria. The use of the proximal point approach is motivated by that of Cabot \cite{Cabot}, who studied hierarchical minimization problems, whose special case is the simple bilevel programming problem. However in Cabot \cite{Cabot} the lower level problem is unconstrained in contrast to our constrained one. Moreover, our convergence analysis does not attempt to generalize Cabot's \cite{Cabot} approach but rather uses a mixture of the approaches due to Cabot \cite{Cabot} and Facchinei  et al. \cite{fpang2011}. In Facchinei et al. \cite{fpang2011} the aim has been to develop an algorithm for the problem VI-C HVI$(F,C,H,\phi)$ as formulated in Section 1. They consider $C$ to be a compact convex set and $F, H$ and $\phi$ are Lipschitz continuous over $C$ with additional monotonicity assumptions. If $H=0$, $F=\nabla g$ and $\phi=f$, then VI-C HVI$(\nabla g,C,0,f)$ reduces to (SBP) with $g$ differentiable. However if $g$ is not differentiable then the (SBP) problem does not belong to the class of VI-C HVI$(F,C,H,\phi)$ problems. The interesting feature of our approach is that even when $g$ is non-differentiable and $C$ is unbounded, we are still able to borrow some ideas from Facchinei  et al. \cite{fpang2011} and adapt them to our circumstances by mixing them with some ideas from Cabot \cite{Cabot}. Note that in our approach even if $ g $ is differentiable, $ \nabla g $ no longer needs to be Lipschitz over $ C $. Let us now briefly outline our motivation for the approach.

Consider a sequence of functions
\begin{eqnarray*}
\psi_n = g + \ve_n f,
\end{eqnarray*}
where $\ve_n >0$ is a 
sequence decreasing to zero. 
Let $x_k$ be the $k$-th iterate  and, motivated by the approach in   proximal point algorithms (see for example the monograph by Burachik and Iusem \cite{regina2011}), 
consider the $(k+1)$-st  iteration as
\begin{eqnarray*}
x_{k+1} \in \eta_k\mbox{-argmin}_C \{\psi_k + \frac{1}{2 \lambda_k} \|.-x_k\|^2\},
\end{eqnarray*}
with $\quad \lambda_k>0, \eta_k>0$, where $\eta\mbox{-argmin}_C h$ denotes the set of $\eta$-minimizers of $h(x)$ subject to $x\in C$. The well known necessary and sufficient optimality condition for the existence of $\eta_k$-minimizer shows that there exists scalars $\eta_k^1 \geq 0$, $\eta_k^2 \geq 0$, $\beta_k^1 \geq 0$ with $\eta_k^1 + \eta_k^2 +\beta_k^1 = \eta_k$ such that
\begin{eqnarray}\label{5.16}
0 \in \partial_{\eta_k^1} \psi_k(x_{k+1}) + \partial_{\beta_k^1}(\frac{1}{2\lambda_k} \|.-x_k\|^2)(x_{k+1}) + N_C^{\eta_k^2}(x_{k+1}).
\end{eqnarray}
Now we know that
\begin{eqnarray*}
\nabla(\frac{1}{2\lambda_k} \|.-x_k\|^2)(x_{k+1})= \frac{1}{\lambda_k}(x_{k+1}-x_k) \in \partial_{\beta_k^1}(\frac{1}{2 \lambda_k}\|.-x_k\|^2)(x_{k+1}).
\end{eqnarray*}
The crux of the idea rests in assuming that (\ref{5.16}) holds as 
\begin{eqnarray*}
0 \in \partial_{\eta_k^1} \psi_k(x_{k+1})+ \left\{\frac{1}{\lambda_k}(x_{k+1}-x_k)\right\} + N_C^{\eta_k^2}(x_{k+1}).
\end{eqnarray*}
Hence,
\begin{eqnarray*}
-\frac{1}{\lambda_k}(x_{k+1}-x_k) \in  \partial_{\eta_k^1} \psi_k(x_{k+1}) + N_C^{\eta_k^2}(x_{k+1}).
\end{eqnarray*}
Naturally $\eta_k^1 + \eta_k^2 \leq \eta_k$. This is used as the principal iteration step in our algorithm.

We are now in a position to state our algorithm and the relevant assumptions needed for the convergence.

\medskip
\begin{tcolorbox}
\noindent{\bf Algorithm for the (SBP) problem}

Let $\psi_n:= g + \ve _nf ,\quad \ve_n >0$.

Step 0. Choose $x_0 \in C$, 
$\epsilon_0 \in (0,\infty)$ and let $k:=0$.

Step 1. Given $x_k$, $\lambda_k$ and $\ve_k$, choose $x_{k+1} \in C$ such that
\begin{eqnarray*}
-(\frac{x_{k+1}-x_k}{\lambda_k}) \in \partial_{\eta_k^1}(\psi_k)(x_{k+1}) + N_C^{\eta_k^2}(x_{k+1})
\end{eqnarray*}
where $\eta_k^1, \eta_k^2 \geq 0$ and $\eta_k^1 + \eta_k^2 \leq \eta_k$.
\end{tcolorbox}
In this algorithm we assume the following hypotheses on the non-negative sequences $\{\epsilon_n\}$, $\{\lambda_n\}$ and $\{\eta_n\}$:

$(\textsl{H}_{\varepsilon})$ The sequence $\{\varepsilon_n\}$ is decreasing 
and $\lim \limits_{n \rightarrow \infty} \varepsilon_n = 0$,

$(\textsl{H}_{\lambda})$ There exists $\underline{\lambda}>0$ and $\overline{\lambda}>0$ such that $\underline{\lambda}\leq \lambda_n \leq \overline{\lambda}$ for all $n \in \mathbb{N}$,

$(\textsl{H}_{\eta})$ The sequence $\{\eta_n\}$ is summable, i.e. $\sum\limits_{n=0}^{\infty} \eta_n <+\infty$.

\noindent Since we are seeking to develop an algorithm for the problem (SBP) under more relaxed assumption than Solodov \cite{Solodov1}, it is natural that the convergence analysis will be more complicated. So we present below some Lemmas which will play a crucial role in the convergence analysis. \\
The following lemma which is used in the proof of the main convergence theorem is a slightly modified version of the Lemma 3.3 given by Cabot\cite{Cabot}. We provide the proof for completeness.

\begin{Lemma}\label{bddlemma}
Let $C \subset \mathbb{R}^n$ be a non-empty, closed, convex set. $g: \mathbb{R}^n \rightarrow \mathbb{R}$ and $f: \mathbb{R}^n \rightarrow \mathbb{R}$ are convex functions which are bounded below. Assume that $S_0 := \arg\min\limits_C g$ and that the set $S_1 := \arg\min \limits_{S_0} f$ is nonempty and bounded. Then for every  $(M_0, M_1) \in \mathbb{R}^2$, the set $ K : = \{x\in C : g(x) \leq M_0\} \cap \{x \in C : f(x) \leq M_1\}$ is bounded.
\end{Lemma}

\noindent{\bf Proof.}
By our assumption,
\begin{eqnarray*}
S_1= \{x\in C: g(x) \leq \min \limits _C g\} \cap \{x\in C: f(x) \leq \min \limits _{S_0} f\}.
\end{eqnarray*}
Since the solution set $S_1$ is bounded, by the property of the recession cone we know that $S_1^{\infty}= \{0\}$. Therefore using the well known recession cone property for closed convex sets (see \cite{rocw}), we get that
\begin{eqnarray}\label{5.2}
S_1^{\infty} = \{x \in \R^n: f^{\infty}(x) \leq 0\} \cap \{x \in \R^n: g^{\infty}(x) \leq 0\} \cap C^{\infty}= \{0\}.
\end{eqnarray}
Now, for any $(M_0,M_1) \in \R^2$, the set
\begin{eqnarray*}
K = \{x \in C: g(x) \leq M_0\} \cap \{x \in C: f(x) \leq M_1\}
\end{eqnarray*}
is bounded if and only if
\begin{eqnarray*}
\{x\in \R^n: g^{\infty}(x) \leq 0\} \cap \{x \in \R^n: f^{\infty}(x) \leq 0\} \cap C^{\infty} = \{0\}.
\end{eqnarray*}
Therefore, from (\ref{5.2}) we can say that the set K is bounded for any $(M_0,M_1) \in \R^2$. \qed

\begin{Lemma}\label{limex}
Let $ \{ h_k \}$ be a sequence of non-negative real numbers. Assume that there exists $ n_0 \in \mathbb{N} $, $ \lambda >0 $ such that
\begin{eqnarray}
h_{k+1} - h_k \leq \lambda \eta_k \quad \forall k \geq n_0,
\end{eqnarray}
where $ \{ \eta_k \} $ is a non-negative summable sequence. Then $ \lim\limits_{ k \rightarrow \infty } h_k $ exists.
\end{Lemma}

\proof
By our assumption, $ h_{k+1} - h_k \leq \lambda \eta_k $ for all $ k \geq n_0 $. Setting $ k = n $ to $ k = n+p-1 $ and then adding up the inequalities for any $ p \in \mathbb{N} $ and $ n \geq n_0 $ we get,
\begin{eqnarray}\label{eqn1}
h_{n+p}-h_n \leq  \lambda\sum\limits_{k=n}^{n+p-1} \eta_k \leq \lambda \sum\limits_{k=n}^{\infty} \eta_k.
\end{eqnarray}
As $ \{ \eta_k \} $ is a summable sequence; for any $ \varepsilon >0 $, there exists $N \in \N $ $( N \geq n_0 )$ such that
\begin{eqnarray*}
\sum\limits_{k=n}^{\infty} \eta_k < \ve,  \quad \forall  n \geq N.
\end{eqnarray*}
Then from (\ref{eqn1}) we get that for all $ n \geq N $,
\begin{eqnarray}\label{0.1}
0 \leq h_{n+p} \leq h_n + \lambda\ve,  \quad \forall  p \in \N.
\end{eqnarray}
Hence $\{h_n\}$ is a bounded sequence. Now we will show that all the subsequential limits of $\{h_n\}$ are the same, leading to existence of
$\lim\limits_{n \rar \infty} h_n$.

\noindent Let us assume on the contrary that there exist two subsequences $ \{h_{n_k}\} $ and $ \{h_{m_k}\} $ of $ \{h_n\} $ converging to different limit points, say $ l $ and $ m $ respectively. Without loss of generality, we   assume that $ l > m $.
Take $\ve= \frac{l-m}{2(\lambda+2)}$, then there exist $N_1 \in \N$ and $M_1 \in \N$ such that
\begin{eqnarray*}
l-\ve &\leq h_{n_k} &\leq l+\ve,  \quad \forall k \geq N_1,
\end{eqnarray*}
and
\begin{eqnarray*}
m-\ve &\leq h_{m_k} &\leq m+\ve, \quad \forall  k \geq M_1.
\end{eqnarray*}
Now we can choose $n_p$ and $m_q$ such that $n_p >m_q \geq N$ and $p,q \geq \max\{N_1,M_1\}$. Then from (\ref{0.1}) we get
\[
h_{n_p} \leq h_{m_q} + \lambda \ve ,
\]
implying that
\[
l-\ve \leq h_{n_p} \leq h_{m_q} + \lambda\ve \leq m+ (\lambda +1)\ve.
\]
As  $\ve=\frac{l-m}{2(\lambda+2)}$ from the last inequality we get that $ m \geq l$, leading to a contradiction. Consequently, $\lim_{k\to\infty} h_k$ exists.\qed

\begin{Lemma}\label{limz}
Let $ I $ be an infinite subset of $ \mathbb{N} $ and $ \{ h_k: k \in I \} $ is a sequence of non-negative real numbers. Assume further that there exists another infinite set $ J \subset I $ such that $ \lim \limits_{k \rightarrow \infty, k \in J } h_k =0 $. Also assume that for any $ k \notin J $, we have
\begin{eqnarray}\label{con}
h_k - h_{k-1} \leq \lambda \eta_{k-1};
\end{eqnarray}
where $ \lambda >0 $ and $ \{\eta_k\} $ is a non-negative, summable sequence. Then $ \lim \limits_{k \rightarrow \infty, k \in I } h_k =0 $.
\end{Lemma}

\proof
Let $\varepsilon>0$ be fixed. Then by our assumptions there exists $ N_\varepsilon \in J $ such that
\begin{eqnarray*}
h_k \leq \frac{\varepsilon}{2} \quad \forall k \geq N_{\varepsilon} \quad \mbox{and} \quad k \in J
\end{eqnarray*}
and
\begin{eqnarray*}
\sum\limits_{k=m}^{\infty} \eta_k \leq \frac{\varepsilon}{2 \lambda} \quad \forall m \geq N_{\varepsilon}.
\end{eqnarray*}
Now for any $ n \in I $ $(n \geq N_{\varepsilon})$, if $ n \notin J $ there exists $ k_n\in J$ such that $k_n <n$ (take the largest of such $k_n$). Since $ J $ is an infinite subset of $ I $, $ k_n \rightarrow \infty $ as $ n \rightarrow \infty $. Then using (\ref{con}) we get that
\begin{eqnarray}\label{x}
0\leq h_n\leq h_{k_n}+ \lambda \sum\limits_{j=k_n}^n \eta_j.
\end{eqnarray}
Note that if $ n \in J $, then the above inequality holds true for $ n = k_n $. Also note that, for all $ n \geq N_{ \varepsilon } $, by the construction of $ k_n $, we can say that $ k_n \geq N_{ \varepsilon }$.
Then for any $ n \in I $ such that $ n \geq N_{ \varepsilon } $, from (\ref{x}) we get
\begin{eqnarray*}
0 \leq h_n \leq \varepsilon.
\end{eqnarray*}
Therefore, $ \lim \limits_{k \rightarrow \infty, k \in I } h_k =0 $.\qed

\begin{Theorem}\label{mainthm1}
Let $ C\subset \mathbb{R}^n $ be a non-empty closed convex set, $g : \mathbb{R}^n \rightarrow \mathbb{R}$ 
and $f : \mathbb{R}^n \rightarrow \mathbb{R}$ be convex functions. Assume that the functions $g$ and $f$ are bounded below, that the set $S_0 := \arg\min\limits_C g$ is nonempty, and that the set $S_1 := \arg\min \limits_{S_0} f$ is nonempty and bounded. Let   $\{ \varepsilon_n\}$  be a sequence satisfying $(\textsl{H}_{\varepsilon})$ and $\sum \limits_{n=0}^{\infty} \varepsilon_n = + \infty$. Further $\{\lambda_n\}$, $\{\eta_n\}$  be   nonnegative sequences   verifying respectively,  $(\textsl{H}_{\lambda})$ and $(\textsl{H}_{\eta})$. Then any sequence $\{x_n\}$ generated by the algorithm satisfies
\begin{eqnarray*}
\lim \limits_{n \rightarrow \infty} d(x_n, S_1)= 0.
\end{eqnarray*}
\end{Theorem}

\proof
We have
\begin{eqnarray*}
\frac{x_{k}-x_{k+1}}{\lambda_k} \in \partial_{\eta_k^1}(g+\varepsilon_k f)(x_{k+1}) + N_C^{\eta_k^2}(x_{k+1}).
\end{eqnarray*}
This implies that there exists $\xi_{k+1} \in N_C^{\eta_k^2}(x_{k+1})$ such that
\begin{eqnarray*}
\frac{x_k-x_{k+1}}{\lambda_k}- \xi_{k+1} \in \partial_{\eta_k^1}(g + \varepsilon_k f)(x_{k+1}).
\end{eqnarray*}
Then for all $x \in C$, by the definition of an $ \eta_k^1 $-subdifferential we have
\begin{eqnarray*}
(g+ \varepsilon_k f)(x) \geq (g+\varepsilon_k f)(x_{k+1}) + \bigg\langle \frac{x_k- x_{k+1}}{\lambda_k}- \xi_{k+1}, x-x_{k+1}\bigg  \rangle - \eta_k^1.
\end{eqnarray*}
Since $\xi_{k+1} \in N_C^{\eta_k^2}(x_{k+1})$, we get $\langle -\xi_{k+1} ,x- x_{k+1}\rangle \geq - \eta_k^2$ for all $x \in C$, which implies that
\begin{eqnarray}\label{a}
(g+ \varepsilon_k f)(x) \geq (g+\varepsilon_k f)(x_{k+1}) + \bigg \langle \frac{x_k- x_{k+1}}{\lambda_k}, x-x_{k+1} \bigg \rangle-\eta_k,
\end{eqnarray}
as $ \eta_k^1 + \eta_k^2 \leq \eta_k $.
In particular for $x= P_{S_1}(x_k)$, we have $ g(x) = \min\limits_C g $ and $ f(x)= \min \limits _{S_0} f $. Hence putting $ x= P_{S_1}(x_k) $ in (\ref{a}), we have
\begin{eqnarray}\label{7.15}
\frac{1}{\lambda_k} \langle x_{k+1}-x_k, x_{k+1}- P_{s_1}(x_k) \rangle \leq \min\limits_C g -g(x_{k+1}) +\nonumber\\
 \varepsilon_k [\min \limits _{S_0} f- f(x_{k+1})] +\eta_k.
\end{eqnarray}
Taking $h_k= \frac{1}{2} d(x_k, S_1)^2$, we get that
\begin{eqnarray*}
 h_{k+1}- h_k &=& \frac{1}{2}\|x_{k+1}- P_{S_1}(x_{k+1})\|^2 - \frac{1}{2}\|x_k- P_{S_1}(x_k)\|^2\\
&\leq & \frac{1}{2}\|x_{k+1}- P_{S_1}(x_k)\|^2 - \frac{1}{2}\|x_k- P_{S_1}(x_k)\|^2\\
&\leq&-\frac{1}{2}\|x_{k+1}- x_k\|^2 + \langle x_{k+1}-x_k, x_{k+1}-P_{S_1}(x_k) \rangle .
\end{eqnarray*}
Then (\ref{7.15}) implies  that,
\begin{eqnarray}\label{7.16}
h_{k+1}-h_k \leq -\frac{1}{2}\|x_{k+1}-x_k\|^2 + \lambda_k [\min\limits_{C}g-g(x_{k+1})] + \\ \nonumber
\lambda_k\varepsilon_k[\min\limits_{S_0}f- f(x_{k+1})] +\eta_k\lambda_k.
\end{eqnarray}
As $x_{k+1} \in C$, $ \min\limits_{C}g \leq g(x_{k+1})$. Therefore by ($H_\lambda$),
 \begin{eqnarray}\label{7.17}
h_{k+1}- h_k \leq -\frac{1}{2}\|x_{k+1}- x_k\|^2 + \lambda_k \varepsilon_k (\min\limits_{S_0}f - f(x_{k+1})) + \overline{\lambda} \eta_k.
\end{eqnarray}

\noindent  Following Cabot \cite{Cabot} we consider the following cases.

\noindent  (a) There exists $k_0 \in \mathbb{N}$ such that for all $ k \geq k_0 $, \
$\min\limits_{S_0} f \leq f(x_k)$,  \text{ and }

\noindent  (b) For all $ k_0 \in \mathbb{N}$ there exists $k \geq k_0$ such that
$f(x_k) < \min\limits_{S_0} f$.

\noindent  {\bf  Case (a)}: We get  from (\ref{7.17}) that  for all $k \geq k_0$,
\begin{eqnarray}\label{case a}
h_{k+1}- h_k \leq \bar{\lambda}\eta_k.
\end{eqnarray}
Then by Lemma \ref{limex}, $ \lim \limits_{k \rightarrow \infty } h_k $ exists and hence the sequence $ \{ h_k \} $ is bounded. Note that, as $ h_k = \frac{1}{2} \| x_k-P_{S_1}(x_k) \|^2 $, the sequence $\{x_k-P_{S_1}(x_k)\}$ is also bounded . Again as $S_1$ is bounded, we have $\{P_{S_1}(x_k)\}$ to be a bounded sequence. Also note that
\begin{eqnarray*}
\|x_k\| \leq \|x_k-P_{S_1}(x_k)\| + \|P_{S_1}(x_k)\|,
\end{eqnarray*}
which clearly shows that $ \{x_k\} $ is a bounded sequence. Thus, there exists a subsequence of $ \{x_k\} $ again say denoted by $ \{x_k\} $ converging to $ \bar{x} $ .\\
We now show  that $ \lim\limits_{k \rightarrow \infty} h_k = 0 $. Observe that from (\ref{7.17}) using $(\textsl{H}_{\lambda})$, we get that
\begin{eqnarray*}
h_{k+1}- h_k \leq \underline{\lambda}\epsilon_k (\min\limits_{S_0}f- f(x_{k+1})) +\overline{\lambda} \eta_k \quad \forall k \geq k_0.
\end{eqnarray*}
Take $ n > k_0 $. Then adding these inequalities for $k=k_0$ to $(n-1)$ one gets
\begin{eqnarray*}
h_n - h_{k_0} \leq \underline{\lambda} \sum\limits_{k=k_0}^{n-1} \epsilon_k [\min\limits_{S_0}f - f(x_{k+1})] + \overline{\lambda}\sum\limits_{k=k_0}^{n-1}\eta_k,
\end{eqnarray*}
which yields
\begin{eqnarray}\label{b}
\underline{\lambda} \sum\limits_{k=k_0}^{n-1} \epsilon_k [f(x_{k+1})- \min\limits_{S_0}f] \leq h_{k_0} -h_{n}+ \overline{\lambda} \sum\limits_{k=k_0}^{n-1} \eta_k .
\end{eqnarray}
As $ \{ \eta_k \} $ is a non-negative summable sequence, let $ \sum\limits_{k=0}^{\infty} \eta_k = C_0 $. Then $\sum\limits_{k=k_0}^{n-1} \eta_k < C_0 $. Also we know that $ h_n \geq 0 $ and $ C_0 < + \infty $. Then from (\ref{b}) we can deduce that
\begin{eqnarray*}
\underline{\lambda} \sum\limits_{k=k_0}^{n-1} \epsilon_k [f(x_{k+1})- \min\limits_{S_0}f] \leq h_{k_0} + C_0.
\end{eqnarray*}
Now as $ \underline{\lambda} >0 $ and $ h_{k_0} + C_0 \leq \infty $, letting  $n \rightarrow \infty$,  we get
\begin{eqnarray*}
\sum\limits_{k=k_0}^\infty \epsilon_k [f(x_{k+1})- \min\limits_{S_0}f] < +\infty .
\end{eqnarray*}
Also note that $ f(x_k) \geq \min_{S_0} f $ for all $ k \geq k_0 $ and $ \varepsilon_k \geq 0 $ for all $ k \in \mathbb{N} $, hence we can conclude that
\begin{eqnarray}\label{7.18}
0 \leq \sum\limits_{k=k_0}^\infty \epsilon_k [f(x_{k+1})- \min\limits_{S_0}f] < +\infty ,
\end{eqnarray}
which implies that
\begin{eqnarray}\label{c}
\lim_{k\to\infty}\epsilon_k [f(x_{k+1})- \min\limits_{S_0}f]=0.
\end{eqnarray}
As $ \lim\limits_{ k \rightarrow \infty } h_k $ exists, we have $\lim_{k\to\infty} (h_{k+1}-h_k)=0$. Then using \eqref{7.17}, (\ref{c}) and $(H_\eta)$ we get
\begin{eqnarray*}
\lim_{k\to\infty}\|x_{k+1}-x_k\|=0.
\end{eqnarray*}
Now, as  $\min\limits_C g\leq g(x_k)$ and $ \underline{\lambda} \leq \lambda_k $ for all $ k \in \N $, we have
\begin{eqnarray*}
\lambda_k ( \min\limits_{C} g - g( x_k )) \leq \underline{\lambda} ( \min\limits_{C} g - g( x_k )).
\end{eqnarray*}
Using this and the fact that $ \min\limits_{S_0} f \leq f(x_k) $ for all $ k \geq k_0 $, we get that
\begin{eqnarray}\label{d}
h_{k+1}-h_k \leq -\frac{1}{2}\|x_{k+1}-x_k\|^2 + \underline{\lambda} [\min\limits_{C}g-g(x_{k+1})] + \eta_k\lambda_k \ \forall k \geq k_0.
\end{eqnarray}
Note that as $ \lim\limits_{k \rightarrow \infty } \| x_{k+1} - x_k \| =0 $, we have $ \lim\limits_{ k \rightarrow \infty } x_{k+1} = \lim\limits_{ k \rightarrow \infty } x_{k} = \bar{x} $. Therefore as $ k \rightarrow \infty $, we get $ 0 \leq \underline{\lambda} [\min\limits_{C} g - g( \bar{x} )]$, which implies that $ g(\bar{x})  \leq \min\limits_C g$ as $ \underline{\lambda} >0 $. As $ C $ is closed, $ \bar{x} \in C $ and hence $ \bar{x} \in \mbox{argmin}_C g $.\\
If possible now assume that $ f(\bar{x})>\min_{S_0}f $, i.e.,   $ \bar{x} \not\in S_1.$ Then there exists $\alpha>0$ such that $f(x_{k+1})-\min_{S_0}f\geq \alpha$ for all sufficiently large $k$, say for $ k \geq k_1 \geq k_0 $. Now using \eqref{7.18},
\begin{eqnarray*}
\alpha \sum\limits_{k=n_0}^\infty \epsilon_k \leq \sum\limits_{k=n_0}^\infty \epsilon_k [f(x_{k+1})- \min\limits_{S_0}f] < \infty,
\end{eqnarray*}
which violates our assumption that $\sum \limits_{k=0}^{\infty} \epsilon_k = + \infty$. Therefore $ f(\bar{x}) = \min_{S_0}f $, which implies that $ \bar{x} \in S_1 $ and finally leads to the fact that $\lim\limits_{k\to\infty} h_k=0$.

\noindent {\bf  Case (b)}:  Let us assume that for each $k_0 \in \mathbb{N}$ there exists $k \geq k_0$ such that $f(x_k) < \min\limits_{S_0}f$. Now we define a sequence $\{\tau_n\}$ as
\begin{eqnarray*}
\tau_n = \max\{k \in \mathbb{N} \, :\,  k \leq n \ \text{ and } \  f(x_k) < \min\limits_{S_0} f\}.
\end{eqnarray*}
We need to first focus on whether the sequence $\tau_n$ is well-defined. It is important to note that the sequence $\{\tau_n\}$ should be understood in the following way.\\
Let $ n_0 \in \mathbb{N}$ be given, then by our assumption there exists $n \in \mathbb{N} $ $( n \geq n_0 ) $ such that
$f(x_n) < \min\limits_{S_0} f $. Now let $ \hat{n} $ be the smallest of all such $ n \geq n_0$ such that $f(x_n) < \min\limits_{S_0} f$. Then $\tau_n$ is well-defined for all $n \geq \hat{n}$. Let us now define the set
\begin{eqnarray*}
I := \{\tau_n: n \geq \hat{n} \}.
\end{eqnarray*}
It is not difficult to see that by our assumptions $\lim\limits_{n \rightarrow \infty} \tau_n = +\infty$. Suppose now that $ \tau_n \leq n-1 $. From the definition of $\tau_n$, for all $k$ satisfying  $\tau_n \leq k \leq n-1$, we have $f(x_{k+1}) \geq \min\limits_{S_0}f$ and hence,  $h_{k+1}-h_k \leq \bar{\lambda} \eta_k$ from (\ref{7.16}).\\
By adding these $(n-\tau_n)$ inequalities, we obtain
$h_n - h_{\tau_n} \leq \bar{\lambda} \sum\limits_{k= \tau_n}^{n-1} \eta_k $,
or, equivalently,
\begin{eqnarray}\label{7.20}
0 \leq h_n \leq h_{\tau_n} + \bar{\lambda}\sum\limits_{k=\tau_n}^{\infty}\eta_k.
\end{eqnarray}
The last inequality is also true for $\tau_n=n$   as  $\bar{\lambda} \sum\limits_{k= n}^{\infty} \eta_k \geq 0$.

 $\bullet$   If we could prove that $\lim\limits_{n \rightarrow \infty}h_{\tau_n}=0$ then from (\ref{7.20}) we can conclude that $\lim\limits_{n \rightarrow \infty}h_n=0 $. Let us consider a subset of $ I $,
\begin{eqnarray*}
J :=\{ k \in I : -\frac{1}{2}\|x_k-x_{k-1}\|^2 + \lambda_{k-1}\epsilon_{k-1}[\min\limits_{S_0}f- f(x_k)]\geq 0 \}.
\end{eqnarray*}

$ \bullet $ If $ J $ is finite, there exists $ \hat{k} \in \N $ such that for all $ k \geq \hat{k} $, $ k \notin J $. Then we have for all $ k \geq \hat{k} $
\begin{eqnarray}
-\frac{1}{2}\|x_k-x_{k-1}\|^2 + \lambda_{k-1}\epsilon_{k-1}[\min\limits_{S_0}f- f(x_k)] < 0.
\end{eqnarray}
Then from (\ref{7.17}) we get that $ h_k -h_{k-1} \leq \bar{\lambda} \eta_{k-1} $. Therefore by Lemma \ref{limex} we can conclude that$ \lim\limits_{ k \rightarrow \infty } h_k $ exists.
This implies $\lim\limits_{k\rightarrow \infty}|h_{k}-h_{k-1}|=0 $
and, in particular, $\lim \limits_{k \in I, k\to\infty}|h_k-h_{k-1}|=0.$ Hence,  $\lim\limits_{k \in I, k\to\infty}\|x_k-x_{k-1}\|=0$, too. Also note that convergence of the sequence $\{h_k\}$ leads to the fact that it is bounded also. Along with this fact the boundedness of the solution  set implies that $\{x_k\}$ is also bounded. Then $\{x_k\}$ has a convergent subsequence. Let us denote the convergent subsequence again by $\{x_k\}$ which converges to $\bar{x} \in C $.\\
Using $ k\not\in J $ and \eqref{7.16} we derive that for all $ k \geq \hat{k} $
\begin{equation*}\label{abc}
h_{k}-h_{k-1}\leq \lambda_{k-1}[\min\limits_C g-g(x_k)]+\eta_{k-1}\lambda_{k-1} \leq \underline{\lambda}[\min\limits_C g-g(x_k)] +\eta_{k-1}\lambda_{k-1}.
\end{equation*}
Then as $ k \to \infty $, we get $ 0 \leq \min\limits_C g - g( \bar{x} ) $ and hence $ \bar{x} \in S_0$. Since $k \in I$, $f(x_k) \leq \min\limits_{S_0} f$. Then continuity of $f$ implies $f(\bar{x}) \leq \min\limits_{S_0}f$.
Therefore,  $\bar{x}$ belongs to the solution set, which means that every limit point of $\{x_{\tau_n}\}$ is a solution of the mentioned simple bilevel optimization problem. Hence, $\lim\limits_{n \rightarrow \infty } h_{\tau_n}=0$.\\
$ \bullet $ If $ J $ is not finite, for all $ k \in J $, $ f(x_k) \leq \min\limits_{S_0}f $ and $ f $ is assumed to be bounded below say by $ M_0 $, then
\begin{eqnarray}\label{7.21}
0 \leq \min\limits_{S_0}f- f(x_k) \leq \min\limits_{S_0}f-M_0.
\end{eqnarray}
Also for all $k \in J $ we have
\begin{eqnarray*}
\|x_k - x_{k-1}\|^2 \leq 2 \lambda_{k-1}\varepsilon_{k-1} [\min\limits_{S_0}f - f(x_k)].
\end{eqnarray*}
Therefore as $k \rightarrow \infty$  (and $ k \in J $) using (\ref{7.21}) we get
\begin{eqnarray}\label{eqextra1}
\lim\limits_{k \rightarrow \infty}\|x_k - x_{k-1}\|=0,  \quad
\end{eqnarray}
which implies  that
\begin{eqnarray}\label{7.22}
\lim\limits_{k \rightarrow \infty}|h_k - h_{k-1}|=0.
\end{eqnarray}
To see this, let us set
\begin{eqnarray*}
h_k= \frac{1}{2}\|x_k-P_{S_1}(x_k)\|^2=: \phi(x_k).
\end{eqnarray*}
It is clear that $\phi$ is a continuous function  as the projection mapping is a Lipschitz function. Thus given $\epsilon >0$, there exists $\delta>0$ such that
if $\|x_k - x_{k-1}\| < \delta$ holds then one has
\begin{eqnarray*}
\vert \phi(x_k)-\phi(x_{k-1})\vert  \leq \varepsilon.
\end{eqnarray*}
Since  \eqref{eqextra1} holds,   there exists $ K \in \mathbb{N}$ such that
$
\|x_k-x_{k-1}\| \leq \delta$  for all $k \geq K$, and thus
\begin{eqnarray*}
 \vert h_k-h_{k-1}\vert  = \vert \phi(x_k)-\phi(x_{k-1})\vert \leq \varepsilon,  \quad \forall~ k \geq K,
\end{eqnarray*}
yielding \eqref {7.22}, which is desired.

\noindent It follows  from (\ref{7.16}), (\ref{7.22}) and $(\textsl{H}_{\lambda})$ that
\begin{eqnarray}\label{7.23}
\lim \limits_{k \rightarrow \infty}g(x_k)=\min\limits_{C}g,\quad k \in J
\end{eqnarray}
i.e.,  $\{g(x_k)\}_{k \in J}$ is a bounded sequence, say bounded by $M_1$. Then
\begin{eqnarray*}
\{x_k : k \in J \} \subset \{x\in C: f(x)\leq \min\limits_{S_0}f\}\cap \{ x\in C: g(x) \leq M_1\}.
\end{eqnarray*}
Lemma \ref{bddlemma} implies that  $\{x_k\}_{ k \in J }$ is a bounded sequence and hence having a convergent subsequence denoted by $\{x_{k_j}\}$, which converges to $ \widetilde{x} $ (say).
Now by continuity of $ g $ and from (\ref{7.23}), we can say that $ \widetilde{x} \in \arg\min\limits_C g $. Since for all $ k \in J $, $ f(x_k) \leq \min\limits_{S_0} f $, from continuity of $ f $ we get that $ f(\widetilde{x})\leq \min\limits_{S_0} f $. Thus $ \widetilde{x} \in S_1 $, which means that $ \lim\limits_{j \rightarrow \infty}h_{k_j}=0 $. This is true for any subsequence of $ \{x_k\}_{k \in \bar{J}} $. Hence
\begin{eqnarray*}
\lim\limits_{k \rightarrow \infty} h_k=0, \quad k\in J.
\end{eqnarray*}
Now $ J \subset I $ and for all $ k \notin J $, we have $ h_k - h_{k-1} \leq \bar{\lambda} \eta_{k-1} $.
Then by Lemma \ref{limz} we can conclude that $ \lim\limits_{k \rightarrow \infty, k \in I} h_k=0 $.\\
Therefore in both the cases whether the set $ J $ is finite or not we have $ \lim\limits_{k \rightarrow \infty, k \in I} h_k=0 $ and hence, $ \lim\limits_{k \rightarrow \infty } h_k =0 $
as mentioned above and  the proof is complete. \hfill $\Box$

\section{SMPEC, Dual Gap Function and Stopping Criteria}

In this section we shall focus on the reformulation of (SMPEC) as a non-smooth simple bilevel problem (SBP-1) in terms of the dual gap function. We have mentioned in section 1, that it might be difficult to use (SBP-1) in real computation, but can be used for theoretical analysis. To begin with we will demonstrate through the following proposition that the dual gap function $ g_D $ of VI$(F,C)$ can be viewed as a penalty function for the (SMPEC) problem.
\begin{Proposition} \label{proppen}
Consider the problem (SMPEC) where $ C $ is a compact convex set. Let $ g_D $ be the dual gap function associated with VI$(F,C)$. Consider the sequence of problems $ (P_k)$ given as
\begin{eqnarray*}
\min\limits_{x\in C} f(x)+ \mu_k g_D (x),
\end{eqnarray*}
where $ \mu_k > 0 $ and $ \mu_k \rightarrow \infty $ as $ k \rightarrow \infty $. Let $ x_k $ be a solution of $ (P_k)$ for each $ k \in \mathbb{N} $. Then any accumulation point of $ \{ x_k \} $ is a solution of (SMPEC)
\end{Proposition}
\proof Let $ \bar{x} $ be an accumulation point of the sequence $ \{ x_k \} $, which exists as $ \{ x_k \} $ is bounded since $ C $ is compact. Without loss of generality let us assume that $ \lim\limits_{k \rightarrow \infty} x_k = \bar{x} $. Let $ \tilde{x} \in \mbox{argmin}_C g_D $ and hence we have for all $ k \in \mathbb{N} $.
\begin{eqnarray}\label{minc}
f(x_k) + \mu_k g_D(x_k) \leq f(\tilde{x}) + \mu_k g_D(\tilde{x})
\end{eqnarray}
As $ \mbox{argmin}_C g_D = \mbox{Sol(VI}(F,C)) $, we know that $ g_D(\tilde{x}) = 0 $. Thus from (\ref{minc}) we have
\begin{eqnarray}\label{gd}
g_D(x_k) \leq \frac{1}{\mu_k} [ f(\tilde{x}) - f(x_k)].
\end{eqnarray}
Since $ C $ is compact, $ g_D $ is finite and continuous as it is convex. Also since $ f $ is convex on $ \R^n $, $ f $ is continuous and hence as $ k \rightarrow \infty $, $ f(x_k) \rightarrow f(\bar{x}) $ and $ g(x_k) \rightarrow g(\bar{x}) $. Thus from (\ref{gd}), as $ k \rightarrow \infty $ we have
\begin{eqnarray}
g_D(\bar{x}) \leq 0,
\end{eqnarray}
showing that $ g_D(\bar{x}) = 0 $ {\it i.e.} $ \bar{x} \in \mbox{Sol(VI}(F,C)) $ . Now as $ g_D(x_k) \geq 0 $, we have
\begin{eqnarray*}
f(x_k) \leq f(x_k) + \mu_k g_D(x_k).
\end{eqnarray*}
Again using (\ref{minc}) we have $ f(x_k) \leq f(\tilde{x}) $. Hence as $ k \rightarrow \infty $, $ f(\bar{x})\leq f(\tilde{x}) $. Since $ \tilde{x} $ is an arbitrary element of $ \mbox{Sol(VI}(F,C)) $, we conclude that $ \bar{x} $ solves the (SMPEC) problem. \qed \\

\noindent Though $(P_k)$ is a convex problem and even if a black-box computes $ g_D $, it might not be very easy to solve $ (P_k)$ as $ g_D $ is non-smooth. Further this scheme will work well when $ C $ is compact, since it will allow us the computation of $ \partial g_D(x) $, $ x \in \mathbb{R}^n $ more easily. A better way to solve the problem (SBP-1) is to follow the approach taken in section 2. Hence we will consider the penalization scheme using the sequence of functions
\begin{eqnarray*}
\phi_k = g_D + \varepsilon_k f \quad k \in \mathbb{N}.
\end{eqnarray*}
The penalization scheme is same as that of Proposition \ref{proppen} with $ \mu_k $ as $ \frac{1}{\varepsilon_k}$. Now for each $ k \in \mathbb{N} $, we shall minimize $ \phi_k $ over $ C $ by using the proximal point method. Let $ x_k $ be the current iterate. Hence the proximal iteration step is given by
\begin{eqnarray}\label{itstep}
x_{k+1} = \mbox{argmin}_C \{ \phi_k + \frac{1}{2 \lambda_k } \| .-x_k \| ^2\}, \quad \lambda_k > 0
\end{eqnarray}
For the theoretical discussion we consider the exact proximal point approach rather than the inexact one of section 2, where $\eta-$minimizers of $ \phi_k + \frac{1}{2 \lambda_k } \| .-x_k \| ^2 $ was considered. Thus by using the standard optimality conditions for convex optimization, we have
\begin{eqnarray*}
0 \in \partial \phi_k (x_{k+1}) + \frac{1}{\lambda_k} ( x_{k+1} - x_k ) + N_C( x_{k+1})
\end{eqnarray*}
Again using the well known sum rule for subdifferential we have
\begin{eqnarray*}
0 \in \partial g_D (x_{k+1}) + \varepsilon_k \partial f (x_{k+1}) + \frac{1}{\lambda_k} ( x_{k+1} - x_k ) + N_C( x_{k+1})
\end{eqnarray*}
Hence,
\begin{eqnarray}\label{scheme}
-\frac{1}{\lambda_k} ( x_{k+1} - x_k ) \in \partial g_D (x_{k+1}) + \varepsilon_k \partial f (x_{k+1}) + N_C( x_{k+1})
\end{eqnarray}
In fact if there exists $ k \in \mathbb{N} $ such that $ x_{k+1} = x_k $, then we have
\begin{eqnarray*}
0 \in \partial g_D (x_{k+1}) + \varepsilon_k \partial f (x_{k+1}) + N_C( x_{k+1})
\end{eqnarray*}
This shows that
\begin{eqnarray}\label{solcond}
0 \in \partial f (x_{k+1}) + \frac{1}{\varepsilon_k} \partial g_D (x_{k+1}) + N_C( x_{k+1})
\end{eqnarray}
From (\ref{solcond}) we conclude that $ x_{k+1} = x_k $ is indeed the solution of the (SMPEC) problem. However from a more practical point of view it is not a good idea to assume the existence of a $k$, such that $ x_{k+1} = x_k $. We may however make use of the condition (\ref{solcond}) to develop a stopping criteria. From the above discussion, a good stopping criteria is one where $\| x_{k+1} - x_k \| $ is less than a given threshold.\\
In order to carry forward our discussion we will first introduce what we mean by a Lagrange multiplier for the (SMPEC) problem. Let us recall from section 3.2 of part-I \cite{DDDP} that (SBP-1) can be reformulated as the following single level problem, called as (r-SMPEC).
\begin{eqnarray*}
&& \min f(x)\\
&& g_D(x) \leq 0\\
&& x\in C.
\end{eqnarray*}
It is clear that $ x^* $ solves (SMPEC) if and only if $ x^* $ is also a minimizer of (r-SMPEC). We say that $ \lambda \geq 0 $, is a Lagrange multiplier of (SMPEC) at a feasible point $ \bar{x} $ if $ \lambda \geq 0 $, is a Lagrange multiplier of (r-SMPEC) at $ \bar{x}$, {\it i.e.}
\begin{eqnarray}\label{lagmul}
0 \in \partial f(\bar{x}) + \lambda \partial g_D(\bar{x}) + N_C (\bar{x}).
\end{eqnarray}
Note that if $ \bar{x} $ is feasible to the (SMPEC) problem, it is also feasible to (r-SMPEC). We say that a feasible point $ \bar{x} $ of (SMPEC) satisfies the Lagrange multiplier rule if there exists $ \lambda \geq 0 $ such that (\ref{lagmul}) holds. Thus in (\ref{solcond}) we can say that Lagrangian multiplier role holds for $ x_{k+1} $ with $\lambda = \frac{1}{\varepsilon_k} $. Again keeping an eye on the practical perspective we must note that it might not always be possible to find a iterate where the Lagrange multiplier rule (\ref{lagmul}) holds. However one may relax the Lagrangian multiplier rule itself by asking only an approximate version of it to be satisfied.\\
We say that $ \bar{x} $ satisfies an $ \ve-$Lagrangian multiplier rule, with $ \ve>0$, if there exists $ u \in \partial f(\bar{x})$, $ w \in \partial g_D(\bar{x})$ and $ v \in N_C(\bar{x})$, along with $ \lambda \geq 0 $ such that
\begin{eqnarray*}
\| u + \lambda w+ v \| \leq \ve.
\end{eqnarray*}
This will now allow us to introduce a stopping criteria for any algorithm associated with the (SMPEC) problem.

\begin{tcolorbox}
\textbf{Stopping Criteria:} Given a threshold $ \ve_0 >0 $, we shall accept the iterate $ x_k $ as an approximate solution to the (SMPEC) problem if $ x_k $ satisfies the $ \ve_0-$ Lagrange multiplier rule.
\end{tcolorbox}
We shall now show that if for $ k $ sufficiently large, the distance between two consecutive iterates $ x_k $ and $ x_{k+1} $ are bounded by a real number depending on $ k $, then $ x_{k+1} $ satisfies the $ \ve_0-$ Lagrange multiplier rule for a given a threshold $ \ve_0 >0 $. We shall demonstrate this through the following proposition.
\begin{Proposition}
Let us consider the (SMPEC) problem with $ C $, a closed convex compact set. Let $ \lambda_k >0 $ be as in (\ref{itstep}). Further assume that $ \{x_k \} $ be a sequence generated by (\ref{itstep}). Let $\ve_0 $ be the given threshold. If for $ k $ sufficiently large we have
\begin{eqnarray*}
\| x_{k+1} - x_k \| \leq \lambda_k \ve_k \ve_0,
\end{eqnarray*}
then the stopping criteria holds at $ x_{k+1} $, with threshold $ \ve_0$.
\end{Proposition}
\proof We would like to recall that the sequence $ \ve_k >0 $ comes from the expression $ \varphi_k = g_D + \varepsilon_k f $. Thus $ \ve_k \downarrow 0 $ as $ k \rightarrow \infty $. As $ \{x_k \} $ is the sequence generated by (\ref{itstep}), we can conclude that (\ref{scheme}) holds, {\it i.e.}
\begin{eqnarray*}
-\frac{1}{\lambda_k} ( x_{k+1} - x_k ) \in \partial g_D (x_{k+1}) + \varepsilon_k \partial f (x_{k+1}) + N_C( x_{k+1}).
\end{eqnarray*}
Hence there exists $ u_{k+1} \in \partial f(x_{k+1})$, $ w_{k+1} \in \partial g_D(x_{k+1})$ and $ v_{k+1} \in N_C(x_{k+1}) $ such that
\begin{eqnarray*}
-\frac{1}{\lambda_k} ( x_{k+1} - x_k ) = w_{k+1} + \varepsilon_k u_{k+1}+ v_{k+1}
\end{eqnarray*}
Dividing both sides by $ \ve_k > 0 $, and taking the usual Euclidean norm, we obtain
\begin{eqnarray}\label{D}
\frac{1}{\lambda_k \varepsilon_k} \| x_{k+1} - x_k \| = \|\frac{1}{\varepsilon_k} w_{k+1} + u_{k+1}+ \frac{1}{\varepsilon_k} v_{k+1} \|
\end{eqnarray}
Note that $ \frac{1}{\varepsilon_k} v_{k+1} \in N_C(x_{k+1}) $. Thus setting $ \frac{1}{\varepsilon_k} v_{k+1}= v_{k+1}^{\prime} $, we have from (\ref{D})
\begin{eqnarray}\label{E}
\|\frac{1}{\varepsilon_k} w_{k+1} + u_{k+1}+ v_{k+1}^{\prime} \|=\frac{1}{\lambda_k \varepsilon_k} \| x_{k+1} - x_k \|
\end{eqnarray}
Now assume that $ k \in \mathbb{N} $ be such that
\begin{eqnarray*}
\| x_{k+1} - x_k \| \leq \lambda_k \ve_k \ve_0.
\end{eqnarray*}
Then from (\ref{E}), we have
\begin{eqnarray*}
\|\frac{1}{\varepsilon_k} w_{k+1} + u_{k+1}+ v_{k+1}^{\prime} \| \leq \ve_0.
\end{eqnarray*}
This shows that $ x_{k+1} $ satisfies the $ \ve_0-$Lagrange multiplier rule with $ \lambda = \frac{1}{\ve_k}$. Thus the stopping criteria holds at $ x_{k+1} $, with threshold $ \ve_0 $. This completes the proof. \qed\\

\section{ Algorithm for the simple MPEC problem}
In this subsection we would like to derive an algorithm for the simple MPEC problem. Of course the simple MPEC problem (SMPEC) is a special case of the variational inequality constrained hemivariational inequality problem with $ H = 0 $ and $ \varphi = f $. However the algorithm we present here differs significantly from that of Facchinei et al.~\cite{fpang2011}. In their convergence proof Facchinei et al.~\cite{fpang2011} show that  any accumulation point of the iterates is a solution of the problem while we show that
the distance of the iterates from the solution set converges to zero. Further let us note that unlike Facchinei et al.~\cite{fpang2011} we consider the function $f$ to be non-smooth. Also note that our motivation for designing our algorithm comes from the approach that we have used for the simple
bilevel programming problem.

\medskip
\begin{tcolorbox}
\noindent{\bf Algorithm for the simple MPEC problem}

\medskip

Step 0. Choose $x_0 \in C$, $e_0$ and $\varepsilon_0 \in (0,\infty)$ and let $k:=0$,

Step 1. Given $x_k$, choose $x_{k+1} \in C$ such that
\begin{eqnarray*}
-(\frac{x_{k+1}-x_k}{ \lambda_k}) \in F(x_{k+1}) + \varepsilon_k \partial_{\eta_k^1}f(x_{k+1}) + N_C^{\eta_k^2}(x_{k+1})
\end{eqnarray*}
where $\eta_k^1, \eta_k^2 \geq 0$ and $\eta_k^1 + \eta_k^2 \leq \eta_k$.\\
\end{tcolorbox}

\noindent In this algorithm, we will assume the following hypotheses on the non-negative sequences $\{\varepsilon_n\}$, $\{\lambda_n\}$ and $\{\eta_n\}$:\\
$(\textsl{H}_{\varepsilon})$ The sequence $\{\varepsilon_n\}$ is non-increasing and $\lim \limits_{n \rightarrow \infty} \varepsilon_n = 0$,\\
$(\textsl{H}_{\lambda})$ There exists $\underline{\lambda}>0$ and $\overline{\lambda}>0$ such that $\underline{\lambda}\leq \lambda_n \leq \overline{\lambda}$ for all $n \in \mathbb{N}$,\\
$(\textsl{H}_{\eta})$ The sequence $\{\eta_n\}$ is summable i.e. $\sum\limits_{n=0}^{\infty} \eta_n <+\infty$.

\begin{Theorem}\label{mainthm2}
Let $F : \mathbb{R}^n \rightarrow \mathbb{R}^n$ be a continuous and monotone plus map and $f : \mathbb{R}^n \rightarrow \mathbb{R}$ is a coercive convex function. Assume that the function $f$ is bounded below, and  that the solution set $S_0$ of the VI(F,C) is nonempty. Then the set $S_1 := \arg\min \limits_{S_0} f$ is nonempty and bounded.
Furthermore, if     $\{ \varepsilon_n\}$ is a sequence  satisfying $(\textsl{H}_{\varepsilon})$ and $\sum \limits_{n=0}^{\infty} \varepsilon_n = + \infty$, and if    $\{\lambda_n\}$, $\{\eta_n\}$  are   nonnegative sequences    verifying respectively $(\textsl{H}_{\lambda})$ and $(\textsl{H}_{\eta})$ then, any sequence $\{x_n\}$ generated by the algorithm satisfies
\begin{eqnarray*}
\lim \limits_{n \rightarrow +\infty} d(x_n, S_1)= 0.
\end{eqnarray*}
\end{Theorem}

\noindent {\bf Proof.}  We have
\begin{eqnarray*}
-(\frac{x_{k+1}-x_k}{ \lambda_k}) \in F(x_{k+1}) + \varepsilon_k \partial_{\eta_k^1}f(x_{k+1}) + N_C^{\eta_k^2}(x_{k+1}).
\end{eqnarray*}

\noindent This implies that there exists $\xi_{k+1} \in N_C^{\eta_k^2}(x_{k+1}) $ such that
\begin{eqnarray*}
-(\frac{x_{k+1}-x_k}{ \lambda_k}) - F(x_{k+1}) - \xi_{k+1} \in \varepsilon_k \partial_{\eta_k^1}f(x_{k+1}).
\end{eqnarray*}

\noindent Then for all $x\in C$,
\begin{eqnarray*}
\varepsilon_k f(x) \geq \varepsilon_k f(x_{k+1}) + \bigg \langle \frac{x_k-x_{k+1}}{\lambda_k} , x-x_{k+1} \bigg \rangle + \langle F(x_{k+1}), x_{k+1}-x \rangle - \eta_k.
\end{eqnarray*}

\noindent In particular for $ x=P_{S_1}(x_k) $ we get,
\begin{eqnarray}
\frac{1}{\lambda_k} \langle x_{k+1}-x_k , x_{k+1}-P_{S_1}(x_k) \rangle &&\leq \varepsilon_k [\min_{S_0}f - f(x_{k+1})] +\nonumber \\ && \langle F(x_{k+1}), P_{S_1}(x_k)- x_{k+1} \rangle + \eta_k. \label{8.25}
\end{eqnarray}

\noindent Now let $h_k= \frac{1}{2} d (x_k,S_1)^2 $. Then
\begin{eqnarray*}
h_{k+1}-h_k \leq -\frac{1}{2}\|x_{k+1}-x_k\|^2 +\langle x_{k+1}-x_k , x_{k+1} - P_{S_1}(x_k) \rangle ,
\end{eqnarray*}
and it now follows from  (\ref{8.25}) that
\begin{eqnarray}\label{8.26}
h_{k+1}-h_k \leq -\frac{1}{2}\|x_{k+1}-x_k\|^2 + \lambda_k \varepsilon_k [\min_{S_0}f - f(x_{k+1})] +\\ \nonumber
 \lambda_k \langle F(x_{k+1}),P_{S_1}(x_k) - x_{k+1} \rangle +
\lambda_k\eta_k,
\end{eqnarray}
which yields
\begin{eqnarray}\label{8.27}
h_{k+1}-h_k \leq -\frac{1}{2}\|x_{k+1}-x_k\|^2 + \lambda_k \varepsilon_k [\min_{S_0}f - f(x_{k+1})] + \lambda_k\eta_k
\end{eqnarray}
by feasibility for the (SMPEC).
We now distinguish the following cases: \

\noindent (a) There exists $k_0 \in \mathbb{N}$ such that for all $k \geq k_0$, \
$\min\limits_{S_0} f \leq f(x_k)$,

\noindent  (b) For all $k_0 \in \mathbb{N}$ there exists $ k \geq k_0$ such that
$f(x_k) < \min\limits_{S_0} f$.

\noindent{\bf Case (a):}  In this case from (\ref{8.27}) we get,  for all $k \geq k_0$,
\begin{eqnarray*}
h_{k+1}- h_k \leq \bar{\lambda}\eta_k.
\end{eqnarray*}
Using the same ideas as in the proof of Theorem \ref{mainthm1} we can show that $\lim\limits_{k\to\infty} h_k$ exists and the sequence of iteration points $\{x_k\}$ is bounded.

\noindent Consider a subsequence of $\{x_k\}$, still denoted by $\{x_k\}$
such that $\lim\limits_{k \rightarrow \infty} x_{k}= \bar{x}$ (say), where $\bar{x} \in C $.

\noindent Now we have to show that $\bar{x}\in \mbox{sol}(VI((F,C)) = S_0$.
Since $\lim\limits _{k\rightarrow \infty}h_k$ exists, we have $\lim\limits _{ k \rightarrow \infty}( h_{k+1}-h_k)=0$. Then for this case, (\ref{8.27}) implies that $\lim\limits _{k \rightarrow \infty}\| x_{k+1}-x_k\|=0$.

\noindent  We now get from (\ref{8.26}), for all $ k \geq k_0 $,
\begin{eqnarray}\label{8.29}
-\frac{1}{\underline{\lambda}}| h_{ k }-h_{k-1}| \leq \langle F(x_{ k }), P_{S_1}(x_{ k -1})- x_{ k } \rangle + \eta_{ k-1}.
\end{eqnarray}
Let $ z_{ k-1}= P_{S_1}(x_{k-1}) $. Then $\{z_{k-1}\}$ is a bounded sequence as $S_1$ is bounded. Hence,  it has a convergent subsequence, still denoted by $ \{z_{k-1}\} $, converging to $\bar{z}$.
Then from (\ref{8.29}), as k goes to infinity, we get
\begin{eqnarray*}
0 \leq \langle F(\bar{x}), \bar{z}- \bar{x} \rangle.
\end{eqnarray*}
As $\bar{z} \in S_1$,  $ \bar z \in \mbox{sol}(VI(F,C))$,
$0 \leq \langle F(\bar{x}),\bar{x}- \bar{z} \rangle$,
implies that
$\langle F(\bar{x}),\bar{x}- \bar{z} \rangle =0$,
and
as F is monotone plus, one has $F(\bar{x})= F(\bar{z})$.

\noindent  Now, for any $x\in C$,
\begin{eqnarray*}
\langle F(\bar{x}), x- \bar{x} \rangle = \langle F(\bar{x}), x- \bar{z} \rangle + \langle F(\bar{x}), \bar{z} -\bar{x} \rangle
 = \langle F(\bar{z}), x- \bar{z} \rangle \geq 0.
\end{eqnarray*}

\noindent Hence, $\bar{x} \in \mbox{sol}(VI(F,C))=S_0$. Now from (\ref{8.27}), we can conclude that
\begin{eqnarray*}
h_{k+1}-h_k \leq \lambda_k \varepsilon_k [\min_{S_0}f - f(x_{k+1})] + \lambda_k\eta_k.
\end{eqnarray*}
Also we have $ \min\limits_{S_0} f \leq f(x_k) $ for all $ k \geq k_0 $. Then arguing same as Case(a) in Theorem \ref{mainthm1} we get that $ f(\bar{x}) = \min\limits_{S_0} f $.
Thus, $\bar{x} \in S_1$, leading to the fact that $\lim\limits_{n \rightarrow \infty} h _k =0$.

\noindent  {\bf  Case (b):}  We assume that for every $ k_0 \in \mathbb{N}$ there exists $ k \geq k_0 $ such that $ f(x_k) < \min\limits_{S_0}f $. Let us consider the sequence $\{\tau_n\}$ defined by
\begin{eqnarray*}
\tau_n = \max\{k \in \mathbb{N}, k \leq n ~and~ f(x_k) < \min\limits_{S_0} f\}
\end{eqnarray*}
and $I:= \{\tau_n: n \in \mathbb{N}\}$.

\noindent  Using the same set $ J $ and same argument as in the proof of Theorem \ref{mainthm1} we can prove that $\lim\limits_{n\to\infty} h_n=0$ provided $\lim\limits_{n\to\infty}h_{\tau_n}=0$.\\

$ \bullet $ If the set $ J $ is finite, following Case (a) we derive $ h_k-h_{k-1}\leq \bar{\lambda} \eta_{k-1} $ for all $ k $ starting with some $ \hat{k} $. Then by Lemma \ref{limex}, $ \lim\limits_{k\to\infty} h_k $ exists.
This implies that $ \lim\limits_{k\in I, k\to\infty} |h_k-k_{k-1}|=0.$ Hence, $\lim\limits_{k\in I, k\to\infty} \|x_k-x_{k-1}\|=0.$
\noindent For $k \in I$, $f(x_k) \leq \min\limits_{S_0}f$ and $f$ is assumed to be bounded below by $M_0$, then
\begin{eqnarray}\label{8.35}
\min\limits_{S_0}f - f(x_k) \leq \min\limits_{S_0}f - M_0.
\end{eqnarray}

\noindent As $\{x_{\tau_n}\}_{\tau_n \in I}$ is a bounded sequence, there exists a convergent subsequence,   again denoted by $\{x_{\tau_n}\}$,  which converges to $\bar{x}$.
Then from (\ref{8.26}) and (\ref{8.35}), as $k \rightarrow \infty$ we have
\begin{eqnarray*}
0 \leq \lim\limits_{k \rightarrow \infty} \langle F(x_k), P_{S_1}(x_{k-1})- x_k \rangle.
\end{eqnarray*}
which implies that $ \bar{x} \in \mbox{sol}(VI(F,C))= S_0 $, as in Case (a).

\noindent  We now  get  $f(\bar{x}) \leq \min\limits_{S_0}f$  from the fact that   $f(x_{\tau_n}) \leq \min\limits_{S_0}f$ and the  continuity of $f$. \
Therefore,  $\bar{x}$ belongs to the solution set, which means that every limit point of $\{x_{\tau_n}\}$ is a solution of the mentioned SMPEC problem. Hence $\lim\limits_{n \rightarrow \infty } h_{\tau_n}=0$

$ \bullet $ If $ J $ is infinite, following the proof of Theorem \ref{mainthm1} we get that
\[\lim\limits_{k\to\infty, k \in J}\|x_k-x_{k-1}\|=0,\]
\begin{equation}\label{8.32}
\lim\limits_{k\to\infty, k \in J} |h_k-h_{k-1}|=0.
\end{equation}
Now as $ f $ is coercive, for all $ k \in J $, $\{x_k\}$ is a bounded sequence and hence has a convergent subsequence. Let $\{x_{k_j}\}$ be the subsequence converging to say, $\bar{x}$. Then from (\ref{8.26}) and (\ref{8.32}) we get that
\begin{eqnarray}\label{8.33}
0 \leq \lim\limits_{k \rightarrow \infty}\langle F(x_k), P_{S_1}(x_{k-1}) - x_k \rangle.
\end{eqnarray}

\noindent Arguing like in Case (a), we can say that $\bar{x} \in S_0$. Since for all $k \in J \subset I $, $f(x_k) \leq\min\limits_{S_0}f$, from the continuity of $f$ we get that $f(\bar{x})\leq \min\limits_{S_0}f$. Thus $\bar{x} \in S_1$, which means $\lim\limits_{j \rightarrow \infty}h_{k_j}=0$. This is true for any convergent subsequence of $\{x_k\}_{k \in J }$. Hence,
\begin{eqnarray*}
\lim\limits_{k \rightarrow \infty} h_k=0, \quad k\in J.
\end{eqnarray*}
Then using the same argument as in case(b) of Theorem \ref{mainthm1} and Lemma \ref{limz} we can conclude that
\begin{eqnarray*}
\lim\limits_{k \rightarrow \infty} h_k=0, \quad k\in I.
\end{eqnarray*}
Therefore we get that $ \lim\limits_{k \rightarrow \infty }h_k =0$ which completes the proof. \hfill $\Box$
$ $\\

\section{Conclusion}
As we end this article we would like honestly confess that it appeared to us that the general numerical scheme analysed here may be difficult to implement. The difficulty arises because it is difficult to compute the non-smooth tools used in these schemes.\\
Meanwhile during the revision of this paper it was observed in Pandit {\it et al.} \cite{PDR}, that if we consider the lower level objective function to be convex and differentiable then a simple algorithm can be developed based on an Armijo type decrease criteria for which we do not need any Lipschitz gradient assumption. Further in Pandit {\it et al.} \cite{PDR}, numerical experiments have been carried out which demonstrated the efficiency of the algorithm. In Pandit {\it et al.} \cite{PDR} no assumption is made about the differentiability of the upper level convex objective function.\\
The discussion in this article also raises the following question that we will take up as a future work.
\begin{enumerate}
\item It might be interesting to see whether the algorithm developed in section 2 can be modified to make it implementable if we assume the (SBP) problem to have smooth objectives. Note that in that case we can avoid the Lipschitz gradient assumption and also do not use any decrease criteria.
\item It appears now to be natural to see if the algorithm developed in section 4 for the (SMPEC) problem can be modified for the following problem.
\begin{eqnarray*}
 \min f(x), \quad \mbox{subject to } x\in \mbox{Sol(GVI}(T,C))
\end{eqnarray*}
where $f: \R^n \rightarrow \R $ is a convex function and $ T : \R^n \rightrightarrows \R^n $ is a maximal monotone map. For the definition of $ \mbox{GVI}(T,C) $, see Part-I \cite{DDDP}. Our intuition might lead us to think that if $ T $ is para-monotone, then we might be able to develop a convergence analysis for any numerical scheme that we might design for the above problem. Our aim would be to get rid of the para-monotone assumption and hence the monotone plus assumption when $ T $ is single-valued.
\item It is important for us to get rid of the monotone plus assumption in section 4. This will allow us to handle, for example the cases where $ F(x) = Mx +q $ where $ M $ is a skew-symmetric matrix. In this case $ F $ is monotone but not monotone-plus. But such class of problems are important as shown in Example 1.3 in Part-I \cite{DDDP}. However in such a setting the value of $ g_D $ can be computed much easily by minimizing a linear function over the set $ C $, which may even be polyhedral in many circumstances. Thus the alternative approach using the equivalent (SBP) formulation may work in this case.
\end{enumerate}
We believe that the analysis carried out in this paper in a very general setting will allow us to now look at the particular and structured cases for which more flexible and implementable algorithms can be designed.

\end{document}